\documentclass[12pt]{amsart}
\pdfoutput=1
\usepackage{amsmath}
\usepackage{amstext}
\usepackage{amsfonts}
\usepackage{amssymb}
\usepackage{amsthm}
\usepackage{amsrefs}
\usepackage{color}
\usepackage{enumitem}
\usepackage{hyperref}
\usepackage{mathtools}
\usepackage{microtype}
\usepackage{thmtools}
\usepackage{tikz-cd}
\allowdisplaybreaks[3]

\numberwithin{equation}{section}
\mathtoolsset{centercolon}

\definecolor{darkblue}{rgb}{0.0,0.0,0.3}
\hypersetup{colorlinks,breaklinks,linkcolor=red,urlcolor=red,anchorcolor=red,citecolor=red}

\theoremstyle{plain}
\newtheorem{thm}{Theorem}[section]

\newtheorem{cor}[thm]{Corollary}
\newtheorem{prop}[thm]{Proposition}

\theoremstyle{definition}
\newtheorem{defn}[thm]{Definition}
\newtheorem{example}[thm]{Example}

\newtheorem{problem}[thm]{Problem}
\newtheorem{hyperconj}[thm]{Hyperrigidity Conjecture}

\newcommand{\bC}{{\mathbb{C}}}

\newcommand{\bN}{{\mathbb{N}}}

\newcommand{\bR}{{\mathbb{R}}}

\newcommand{\A}{{\mathcal{A}}}
\newcommand{\B}{{\mathcal{B}}}

\newcommand{\I}{{\mathcal{I}}}

\newcommand{\M}{{\mathcal{M}}}
\newcommand{\cM}{{\mathcal{M}}}

\renewcommand{\O}{{\mathcal{O}}}

\renewcommand{\S}{{\mathcal{S}}}

\newcommand{\fA}{{\mathfrak{A}}}
\newcommand{\fB}{{\mathfrak{B}}}

\newcommand{\rA}{{\mathrm{A}}}
\newcommand{\rB}{{\mathrm{B}}}
\newcommand{\rBI}{{\mathrm{B^\infty}}}

\newcommand{\rC}{{\mathrm{C}}}

\newcommand{\fs}{{\mathfrak{s}}}

\renewcommand{\phi}{\varphi}

\newcommand{\Aut}{\operatorname{Aut}}

\newcommand{\ncconv}{\operatorname{ncconv}}

\newcommand{\ip}[1]{\langle #1 \rangle}

\newcommand{\ol}{\overline}

\newcommand{\Ran}{\operatorname{Ran}}

\newcommand{\Rep}{\operatorname{Rep}}
\newcommand{\spn}{\operatorname{span}}

\newcommand{\ucp}{\text{u.c.p.}}

\newenvironment{sbmatrix}{\left[\begin{smallmatrix}}{\end{smallmatrix}\right]}
\newcommand{\sot}{\textsc{sot}}
\newcommand{\wot}{\textsc{wot}}

\newcommand{\qfor}{\quad\text{for}\quad}
\newcommand{\qforal}{\quad\text{for all}\quad}

\newcommand{\AND}{\ \text{and}\ }

\DeclareMathOperator{\graph}{Graph}
\DeclareMathOperator{\epi}{Epi}

\DeclareMathOperator{\ncpmaps}{CP_{nor}}

\DeclareMathOperator{\bor}{Bor}

\newcommand{\env}[1]{\overline{#1}}

\newcommand{\cmax}{\mathrm{C}_{\textup{max}}^*}

\makeatletter
\def\widebreve#1{\mathop{\vbox{\m@th\ialign{##\crcr\noalign{\kern\p@}%
  \brevefill\crcr\noalign{\kern0.1\p@\nointerlineskip}%
  $\hfil\displaystyle{#1}\hfil$\crcr}}}\limits}

\def\brevefill{$\m@th \setbox\z@\hbox{}%
 \hfill\scalebox{1.1}{\rotatebox[origin=c]{90}{(}} \kern4pt $}
\makeatother

\begin{document}

\title[Noncommutative Choquet theory]{Noncommutative Choquet theory:\\A Survey}

\author[K.R. Davidson]{Kenneth R. Davidson}
\address{Department of Pure Mathematics\\ University of Waterloo\\Waterloo, ON, N2L 3G1, Canada}
\email{krdavids@uwaterloo.ca\vspace{-2ex}}

\author[M. Kennedy]{Matthew Kennedy}
\address{Department of Pure Mathematics\\ University of Waterloo\\Waterloo, ON, N2L 3G1, Canada}
\email{matt.kennedy@uwaterloo.ca}

\begin{abstract}
We survey noncommutative Choquet theory and some of its applications.
\end{abstract}

\thanks{First author supported by NSERC Grant Number 2018-03973.}
\thanks{Second author supported by NSERC Grant Number 50503-10787.}
\maketitle
\tableofcontents

The Krein-Milman theorem asserts that a compact convex subset of a locally convex space is generated by its extreme points, in the sense of taking the closed convex hull. On the other hand, Milman's partial converse asserts that the extreme points are an optimal generating set, in the sense that they necessarily belong to the closure of any other generating set. It is therefore natural to ask how arbitrary points in a compact convex set can be obtained from extreme points.

This question is answered by Choquet theory, which provides an integral representation for arbitrary point in a compact convex set by a probability measure concentrated on the extreme points \cite{Ch1956}. More generally, Choquet theory provides a powerful functional-analytic framework for the analysis of convex sets, and has numerous important applications in approximation theory, potential theory, ergodic theory, and many other areas of mathematics.

In the theory of operator algebras, objects often exhibit what can best be described as ``noncommutative'' convex structure. Several attempts have been made to capture this structure formally, most successfully with the theory of matrix convexity. However, despite its success, this theory suffers from a significant issue, which is the lack of an appropriate notion of extreme point.

Recently, the authors introduced a new theory of noncommutative convexity that refines the theory of matrix convexity and, in particular, completely resolves this issue with extreme points. Based on this framework, the authors further developed a corresponding noncommutative Choquet theory, in the process demonstrating that a surprising amount of the classical theory extends to the noncommutative setting.

In this article, we will survey noncommutative convexity and noncommutative Choquet theory, and discuss some applications to operator algebras, group theory and dynamical systems.

Classical Choquet theory exploits the duality between a compact convex set $K$ and the space $\rA(K)$ of continuous complex-valued affine functions on $K$.
This is a unital self-adjoint subspace of the C*-algebra $\rC(K)$ of continuous complex-valued functions on $K$. Since $\rA(K)$ separates points, it generates $\rC(K)$ as a C*-algebra.

A unital self-adjoint generating subspace of a commutative C*-algebra is called a function system. An important classical theorem of Kadison \cite{Kad1951}, often called Kadison's representation theorem, shows that an arbitrary function system  $S$ is unitally order isomorphic to the function system $\rA(K)$ inside $\rC(K)$ as above, where $K$ denotes the state space of $S$, i.e. the set of unital positive linear functionals on $S$ equipped with the relative weak* topology.

The above correspondence between a function system and its state space yields a contravariant categorical equivalence between the category of compact convex sets with continuous affine maps as morphisms, and the category of function systems with unital order homomorphisms as morphisms. 

The continuous convex functions on $K$ play a critical role in Choquet theory. In fact, the interaction between $\rA(K)$, $\rC(K)$ and the continuous convex functions on $K$ provides the analytic machinery from which the convex structure of $K$ can be extracted.

The Choquet order on the probability measures on $K$ is defined by comparing the integrals of continuous convex functions. The probability measures that are maximal in this order turn out to be supported (in a technical sense) on the set of extreme points of $K$, and this yields the Choquet-Bishop-de Leeuw representation theorem \cite{BdL1959}, which asserts that every point in $K$ has an integral representation by a probability measure concentrated on the extreme points.

A particularly important class of compact convex sets are the Choquet simplices, which are characterized by the property that every point has a unique representing probability measure that is maximal in the Choquet order. A Choquet simplex is a generalization of a simplex in finite dimensions.

A Choquet simplex in which the set of extreme points is closed is called a Bauer simplex. The work of Bauer \cite{Bauer} shows that Bauer simplices are precisely the state spaces of unital commutative C*-algebras. A Choquet simplex in which the set of extreme points is dense is called a Poulsen simplex. The work of Lindenstrauss, Olsen and Sternfeld \cite{Lindenstrauss-Olsen-Sternfeld}  shows that there is a unique metrizable Poulsen simplex. Both Bauer simplices and Poulsen simplices play an important role in ergodic theory (see e.g. \cite{Glasner-Weiss}), and the Poulsen simplex is an important example of a Fra\"{i}ss\'{e} limit in mathematical logic (see e.g. \cite{Lupini1}).

In recent years, various notions of ``noncommutative'' convexity have proven important in the theory of operator algebras. The notion of a matrix convex set, which is actually a family of related convex sets in the space of $n \times n$ matrices over a locally convex space, was introduced in the work of Wittstock \cites{Wit1981,Wit1984} on completely positive maps. The theory was further developed by Effros and Winkler \cite{EffWin1997}, who proved an analogue of the Hahn-Banach theorem in this setting, and then by Webster and Winkler \cite{WebWin1999}, who proved an analogue of Kadison's representation theorem.

The theory of noncommutative (nc) convexity was introduced by the authors as a refinement of the theory of matrix convexity. A subtle, but very important difference is that the theory of nc convexity permits the consideration of infinite matrices. One compelling reason for this difference is that there is no adequate notion of extreme point for compact matrix convex sets. Allowing infinite matrices yields the existence, for compact nc convex sets, of a naturally defined and minimal set of extreme points, along with appropriate analogue of the Krein-Milman theorem and Milman's partial converse.

Another motivation for the theory of nc convexity is provided by the theory of nc functions on nc convex sets. The notion of an nc function has a long history, going back to work of Taylor \cite{Tay1972} on multivariable noncommutative functional calculus, the work of Takesaki  \cite{Tak1967} on a noncommutative Gelfand representation for C*-algebras, and the more recent work of Voiculescu \cite{Voi2008} on free analytic functions.

For a compact nc convex set $K$, there is a natural notion of affine nc function on $K$. The space $\rA(K)$, consisting of continuous affine nc functions generates the C*-algebra $\rC(K)$ of continuous nc functions on $K$. The von Neumann algebra $\rB(K)$ of all bounded nc functions on $K$ is isomorphic to the bidual of $\rC(K)$. There is a natural notion of convex nc function and, as in the commutative setting, the interplay between $\rA(K)$, $\rC(K)$ and the continuous convex nc functions on $K$ provides the analytic machinery from which the nc convex structure of $K$ can be extracted.

The theorem of Webster-Winkler \cite{WebWin1999} implies an analogue of Kadison's representation theorem for compact nc convex sets. Specifically, an arbitrary operator system $S$, i.e. a unital self-adjoint subspace of a C*-algebra, is unitally completely order isomorphic to the oeprator system $\rA(K)$ inside $\rC(K)$ as above, where $K$ denotes the nc state space of $S$, i.e. the set of \ucp maps on $S$. This yields a contravariant categorical equivalence between the category of compact nc convex sets and the category of operator systems. 

In this noncommutative setting, measures are replaced by unital completely positive (\ucp) maps on $\rC(K)$ with ranges in $\B(H_n)$ for Hilbert spaces $H_n$ of dimension $n$, where $n$ ranges over cardinal numbers. There is a natural notion of barycentre for such a map, which provides a representation of points in $K$.

There is a natural notion of nc Choquet order on \ucp\ maps on $\rC(K)$, arising from comparison on continuous convex nc functions, so that one can seek maximal maps which represent a point in $K$. We show that this nc Choquet order can equivalently be obtained by considering dilations of \ucp\ maps on $\rC(K)$. This fact has a number of important consequences, and leads to a noncommutative analogue of the Choquet-Bishop-de Leeuw theorem that every point in $K$ can be represented by a map that is supported on the extreme points of $K$ in a certain precise sense. In the separable case, one gets a bona fide integral representation. These results are contained in \cite{DK_ncCh}.

The notion of an nc Choquet simplex was defined and characterized by the second author and Shamovich \cite{Kennedy-Shamovich}. This is a much larger class of sets than the set of classical Choquet simplices, but they satisfy noncommutative analogues of many of the nice properties of classical Choquet simplices. 
An nc Choquet simplex is called a Bauer nc simplex if the set of extreme points is closed.
It is shown that $K$ is a Bauer simplex if and only if $\rA(K)$ is unitally completely order isomorphic to a unital C*-algebra, which is a noncommutative analogue of the classical result of Bauer. There is also a natural noncommutative generalization of a Poulsen simplex \cites{KirWas1998, Lupini2}. 

One of the striking applications of classical Choquet theory is in approximation theory.
Motivated by Bernstein's proof of Weierstrass's theorem for $C[0,1]$, Korovkin observed that any sequence of positive maps from $C[0,1]$ to itself
which converges to the identity map on $\{1,x,x^2\}$ must converge pointwise to the identity on $C[0,1]$.
\u Sa\u skin \cite{Sas1967} established a far-ranging generalization using Choquet theory.
This result is limited to the metrizable case. However using ideas from noncommutative Choquet theory in \cite{DK2021},
we show how to extend these results to compact Hausdorff spaces. 

Arveson \cite{Arv2011} showed that there is a truly non-commutative extension of Korovkin's result.
He formulated the general problem for operator systems in arbitrary separable C*-algebras.
Then he formulated his Hyperrigidity Conjecture, which if true would yield a non-commutative generalization of \u Sa\u skin's theorem.
He showed that this conjecture was valid in a variety of circumstances.
In \cite{DK2021}, we show that for function systems, this comes down to a comparison of the classical and non-commutative Choquet orders on measures on $K$.
Recent work of Bilich and Dor-On \cite{Bilich-DorOn} produced a counterexample to the conjecture.
Clou\^{a}tre and Thompson \cite{Clouatre-Thompson} show that a restricted version of the conjecture does hold.

In the last section, we discuss the application of non-commutative Choquet theory to dynamical systems. One nice application is a a new dynamical characterization of property (T) for locally compact groups \cite{Kennedy-Shamovich} that extends an important result of Glasner and Weiss \cite{Glasner-Weiss}.

There are several other applications of nc convexity and nc Choquet theory that we will only briefly mention here. The categorical equivalence between compact nc convex sets and operator systems, was extended to a duality between (potentially) nonunital operator systems and pointed compact nc convex sets by the second author, Manor and Kim \cite{Kennedy-Manor-Kim}. Humeniuk \cite{Humeniuk} obtained a far-reaching noncommutative generalization of Jensen's inequality. Finally, Humeniuk, Manor and the second author \cite{Humeniuk-Kennedy-Manor} obtained a characterization of dualizability for nonunital operator systems. 

For background in classical Choquet theory, see \cites{Alfsen, Phelps, LMNS}.
For background on C* and von Neumann algebras, see \cites{Davidson, Pedersen, Tak2003}.
For background on operator systems and operator spaces, see \cites{Paulsen, P2003, ER2000}.

\section{Noncommutative convex sets}

We work in an operator space $E$ which is a dual space with predual $E_*$. Then weak-$*$ closed bounded sets are compact.
The space $M_{m,n}(E)$ of $m \times n$ matrices with coefficients in $E$ is a natural construct for an operator space.
One can also talk about infinite matrices $M_{m,n}(E)$ over $E$, for possibly infinite cardinals $m,n$, with a compatible operator space structure \cite{ER2000}*{Section 10.1}.

We fix an infinite cardinal $\kappa$ at least as large as the cardinality of some dense subset of $E$.
Let $\M(E) = \coprod_n M_n(E)$ be the disjoint union over cardinals $n \le \kappa$.
We suppress $\kappa$ in the notation, since it plays no role except to make this a set.
We will write ``for all $n$'' when we mean for all $n \le \kappa$, since it will also follow for larger cardinals if needed.
If $E$ is separable, we can take $\kappa = \aleph_0$. Let $\M = \M(\bC)$.

A subset $K = \coprod_n K_n$ of $\M(E)$ is  \emph{nc convex}  if
\begin{enumerate}
\item $K$ is graded; i.e., $K_n \subset M_n(E)$ for all $n$.
\item If $\{ \alpha_i \in M_{n,n_i} : i \in I \}$ is a family of isometries such that $\sum_I \alpha_i \alpha_i^* = 1_n$
and $\{ x_i \in K_{n_i} : i \in I \}$ is a bounded set, then $x = \sum_I \alpha_i x_i \alpha_i^* \in K_n$.
\item if $\beta \in M_{m,n}$ is an isometry and $x \in K_m$, then $\beta^* x \beta$ belongs to $K_n$.
\end{enumerate}
Say that $K$ is closed or compact if each $K_n$ is closed or compact.

Property (2) says that $K$ is closed under arbitrary direct sums, and (3) says that $K$ is closed under compressions.
If $\{ x_i \in K_{n_i} : i \in \I\}$ is a bounded subset and $\beta_i \in M_{n_i,n}$ such that $\sum _i \beta_i^* \beta_i = 1_n$,  
we can form an \emph{nc convex combination} by $x = \sum_i \beta_i^* x_i \beta_i$.
This belongs to $K_n$ because $y := \bigoplus_i x_i$ belongs to $K_m$ for $m = \sum_i n_i$ by (2), and 
\[
 \beta = \begin{bmatrix} \beta_1 \\ \beta_2 \\ \vdots\end{bmatrix} \in M_{m,n},
 \ \ 
 \beta^* \beta = \sum _i \beta_i^* \beta_i = 1_n
 \ \AND \ 
 \beta^* y \beta =  \sum_i \beta_i^* x_i \beta_i = x.
\] 
That is, $\beta$ is an isometry, so that $x \in K_n$ by (3). Thus an nc convex set is a set which is closed under nc convex combinations. 

The basic example is the nc state space $\S(S)$ of an operator system $S$. 
We let $K_n = UCP(S,M_n) \subset CB(S, M_n) = M_n(S^*)$, where $S^*$ is the operator space dual of $S$ (not the operator system dual).
Each $K_n$ is compact in the point-weak-$*$ topology. It is easy to see that the direct sum of \ucp\ maps and a compression of a \ucp\ map are \ucp\
Thus $\S(S)$ is a compact nc convex set.

A \emph{matrix convex set} $K = \coprod_{n\ge1} K_n$, where $K_n \subset M_n(E)$ for $n < \infty$,
is a set which is graded, closed under finite direct sums and compressions.
A closed nc convex set is uniquely determined by the sets $K_n$ for $n<\infty$.

A function $\theta:K \to L$ between nc convex sets is \emph{nc affine} if it is graded and preserves direct sums and compressions.
It is continuous if it is continuous at each level. Let $\rA(K)$ denote the space of continuous nc affine functions.
If $\theta$ and $\theta^{-1}$ are nc affine, then $\theta$ is an \emph{nc affine homeomorphism}.

The following is taken from Webster-Winkler \cite{WebWin1999} and Davidson-Kennedy \cite{DK_ncCh}.

\begin{thm} 
If $K$ is a compact nc convex set, then $\rA(K)$ is an operator system. The map taking $x\in K$ to $\delta_x$, the evaluation map,
is an nc affine homeomorphism of $K$ onto the nc state space of $\rA(K)$.

Conversely, if $S$ is an operator system with nc state space $K=\S(S)$, then the map from $s \in S$ to the function $\hat s(x) = x(s)$ for $x\in K$
is a complete order isomorphism of $S$ onto $\rA(K)$. 
\end{thm}

If $\theta:K \to L$ is a continuous nc affine map, then we can define a \ucp\ map $A(\theta) : \rA(L) \to \rA(K)$ by $A(\theta)(b)(x) = b(\theta(x))$ for $b \in \rA(L)$ and $x \in K$. 
Conversely, if $\phi:S \to T$ is a \ucp\ map between operator systems and $\S(T) = L$, then
$\S(\phi) (y) (s) = y(\phi(s))$ for $y \in L$ and $s \in S$ is an nc affine homomorphism of $L$ into $K$.

\begin{cor} \label{C:categorical duality}
This provides a contravariant functor $A$ from the category of compact nc convex sets with  morphisms consisting of continuous nc affine maps
onto the category of operator systems with morphisms consisting of \ucp\ maps; and $\S$ is its inverse.
\end{cor}

This identification allows us to apply the notion of dilation to $K$. Say that $y \in K_m$ \emph{dilates} $x\in K_n$ and that  $x$ is a \emph{compression} of $y$ if there is an isometry $\beta \in M_{m,n}$ such that $x = \beta^* y \beta$.
This is a \emph{trivial dilation} if $y \simeq x \oplus z$ with respect to $\Ran \beta \oplus (\Ran \beta)^\perp$.
A point $x$ is \emph{maximal} if it has no non-trivial dilations.

The correspondence between points of $K$ and $\ucp$ maps on $\rA(K)$ connects this notion of maximality with the \emph{unique extension property}.
This means that $x \in K_n$, which corresponds to a \ucp\ map  from $\rA(K)$ to $\B(H_n)$, has a unique completely positive extension to a map from
$C^*(\rA(K))$ and this extension is a representation. We will explain shortly what this C*-algebra is, but the unique extension property is independent of the choice of C*-algebra. 
The Dritschel--McCullough theorem \cite{DriMcC2005} can be restated as: 

\begin{thm} \label {T:DM}
Every $x\in K$ has a maximal dilation $y$.
\end{thm}

The next goal is an analogue of the Krein-Milman theorem. Let $K$ be an nc convex set. 
A point $x \in K_n$ is \emph{nc extreme} if whenever $x$ is written as a finite nc convex combination 
$x = \sum \alpha_i^* x_i \alpha_i$ for $x_i \in K_{n_i}$ and nonzero $\alpha_i \in M_{n_i,n}$ satisfying $\sum \alpha_i^* \alpha_i = 1_n$, 
then each $\alpha_i$ is a positive scalar multiple of an isometry $\beta_i \in M_{n_i,n}$ satisfying $\beta_i^* x_i \beta_i = x$ 
and each $x_i$ decomposes with respect to the range of $\alpha_i$ as a direct sum $x_i = y_i \oplus z_i$ for $y_i,z_i \in K$. 
Note that $y_i$ unitarily equivalent to $x$ via $\beta_i$. The set of all nc extreme points is denoted $\partial K = \coprod_n (\partial K)_n$.

We also define a point $x\in K_n$ to be \emph{pure} if whenever $x$ is written as a finite nc convex combination 
$x = \sum \alpha_i^* x_i \alpha_i$ for $x_i \in K_{n_i}$ and nonzero $\alpha_i \in M_{n_i,n}$ satisfying $\sum \alpha_i^* \alpha_i = 1_n$, 
then each $\alpha_i$ is a positive scalar multiple of an isometry $\beta_i \in M_{n_i,n}$. This is a weakening of the notion of being extreme.
It can be shown that $x$ is extreme if and only if it is pure and maximal.

Again the correspondence between points in $K$ and \ucp\ maps on $\rA(K)$ translates purity to purity of the \ucp\ maps.
So a pure maximal point corresponds to a maximal \ucp\ map with the unique extension property.
Purity means that this unique extension is an irreducible representation. 
These are the \emph{boundary representations} of $\rA(K)$ as defined by Arveson \cite{Arv1969}.
The Davidson--Kennedy theorem \cite{DK2015} can be restated as: 

\begin{thm} \label {T:DK}
Every pure point $x\in K$ has an extreme dilation $y$.
\end{thm}

Now we can state the nc Krein--Milman theorem.

\begin{thm} \label {T:nc krein-milman}
A compact nc convex set is the closed nc convex hull of its nc extreme points; i.e. $K = \ol{\ncconv(\partial K)}$.
\end{thm}

As in the classical case, there is a useful partial converse.

\begin{thm} \label {T:nc milman converse}
Let $K$ be a compact nc convex set, and let $X \subset K$ be a closed set which is closed under compressions.
If\/ $\ol{\ncconv(X)} = K$, then $\partial K \subset X$.
\end{thm}

\section{Noncommutative functions}

Functions on a convex set are replaced by nc functions, which respect some of the structure of nc convex sets.
Let $K$ be a compact nc convex set.
A function $f : K \to \M = \coprod_n M_n$ is an \emph{nc function} if 
\begin{enumerate}
\item $f$ is graded: $f(K_n) \subseteq M_n$ for all $n$,
\item $f$ respects direct sums: if $x_i \in K_{n_i}$ for $i \in I$ and there are isometries $\alpha_i \in M_{n_i,n}$ 
satisfying $\sum_{i \in I} \alpha_i \alpha_i^* = 1_n$, then
\[
  f \big( \sum \alpha_i x_i \alpha_i^* \big) = \sum \alpha_i f(x_i) \alpha_i^* .
\]
\item $f$ is unitarily equivariant: if $x \in K_n$ and $\beta \in M_n$ is unitary, then $f(\beta x \beta^*) = \beta f(x) \beta^*$.
\end{enumerate}
Say that $f$ is \emph{bounded} if $\|f\|_\infty  = \sup_{x \in K} \|f(x)\| < \infty$.
Let $\rB(K)$ denote the space of all bounded nc functions on $K$.

Nc affine functions are nc functions.
However, nc functions do not generally preserve compressions, only unitary equivalence.

\begin{example}
Let $\O_2$ be the Cuntz algebra with generators $\fs_1,\fs_2$, and let $S= \spn\{1, \fs_1,\fs_2,\fs_1^*,\fs_2^*\}$.
Let $K$ be its nc state space of $S$.
Every \ucp\ map $x \in K_n$ sends $\big[\fs_1,\ \fs_2 \big]$ to $T = \big[x(\fs_1),\ x(\fs_2) \big]$, which is a row contraction;
and this row contraction determines $x$. 
Conversely the dilation theorem of Frazho \cite{Fra1982} shows that every $1 \times 2$ row contraction dilates to a row isometry, 
which can then be dilated to a row unitary in various ways, determining a representation of $\O_2$.
So $K$ consists of all $1 \times 2$ row contractions $T$.
There is a unique decomposition $T = U \oplus T'$ with respect to $H_n = M \oplus M^\perp$,
where $U$ is a row unitary and  $T'$ has no row unitary summand.
Define an nc function $f$ by $f(x) = P_M$.
This function is graded and unitarily equivariant.
To check direct sums, note that the direct sum of row contractions with no row unitary summand can have no unitary summand.
The function $f$ is a projection in $\rB(K)$.
A significant feature of this function is that $f|_{K_n} = 0$ for all $n<\infty$. 
So $f$ is not determined by its finite part. 
\end{example}

It is easily seen that we can add, multiply and take the adjoint of bounded nc functions, making $\rB(K)$ into a $*$-algebra.
It comes with a large family of $*$-representations, namely evaluation at points of $K$. Define $\pi_u:\rB(K) \to \prod_n \prod_{x\in K_n} M_n$ by
\[
 \pi_u(f) = \bigoplus_{x\in K} f(x) = f(x_u), \quad\text{where}\quad x_u = \bigoplus_{x\in K} x .
\]
It is straightforward to verify that this representation is isometric and the range is \wot-closed. Hence $\rB(K)$ is a von Neumann algebra.
Define $\rC(K) = C^*(\rA(K))$ as a subalgebra of $\rB(K)$.

\begin{thm} 
$\rB(K)$ is the universal von Neumann algebra $\rC(K)^{**}$ of $\rC(K)$.
The weak-$*$ closure of $\rA(K)$ in $\rB(K)$ equals the space of all bounded nc affine functions on $K$.
\end{thm}

Every operator system $S$ is contained in a universal C*-algebra $\cmax(S)$, meaning that there is a unital completely isometric imbedding $\kappa:S \to \cmax(S) = C^*(\kappa(S))$ with the property that whenever $j:S \to \fB = C^*(j(S))$ is a unital complete isometry, there is a unique surjective $*$-homomorphism $\pi:\cmax(S) \to \fB$ such that $j = \pi\kappa$. It is easily shown that whenever $\phi:S\to \fB$ is a \ucp\ map into a C*-algebra, then there is a unique $*$-homomorphism $\pi_\phi: \cmax(S) \to \fB$ such that $\phi = \pi_\phi \kappa$. This property characterizes $\cmax(S)$. See \cites{Blecher-LeMerdy, KirWas1998}.

Every \ucp\ map from $\rA(K)$ to $M_n$ is evaluation at a point $x \in K_n$. This extends to the normal $*$-homomorphism $\delta_x$ of $\rB(K)$ given by evaluation at $x$, and in particular to a $*$-homomorphism of $\rC(K)$. Every $*$-representation of $\rC(K)$ arises in this way. Thus we get

\begin{cor}
For a compact nc convex set $K$, $\rC(K) \simeq \cmax(\rA(K))$.
\end{cor}

The issue of continuity of functions in $\rC(K)$ is a bit delicate. It turns out that they need not be continuous on $K$ with its usual topology.
We have noted that can identify $K$ with the space $\Rep(\rC(K))$ of $*$-representations of $\rC(K)$.
For any C*-algebra $\fA$, each $a \in \fA$ determines a function $\hat a$ on $\Rep(\fA)$ by $\hat a(\pi) = \pi(a)$ which has a lot in common with nc functions.
In particular, $\hat a$ preserves direct sums and is unitarily equivariant.
This idea is due to Takesaki \cite{Tak1967}, which along with Bichteler \cite{Bic1969}, yields a non-commutative Gelfand representation for an arbitrary C*-algebra.

It can be shown that various natural topologies on this space coincide, including the point-\wot, point-\sot, and point-ultrastrong-$*$ topologies.
The function $\hat a$ is continuous in this topology.
There is a natural map of $\Rep(\rC(K))$ onto $K$ by restriction to $\rA(K)$.
This map is a continuous bijection, but since $\Rep(\rC(K))$ is generally not compact, the inverse is usually not continuous.

Let $\tau_{us*}$ be the weakest topology on $K$ so that every $f\in\rA(K)$ is continuous from $K_n$ into $(M_n, \text{ultrastrong-}*)$ for all $n$. 
Also let $L_1$ denote the classical state space of $\rC(K)$ with its usual weak-$*$ topology.

\begin{thm} \label {T:continuity of C(K)}
Let $f \in \rB(K)$. The following are equivalent:
\begin{enumerate}
\item $f \in \rC(K)$.
\item $\hat f$ is continuous on $(\Rep(\rC(K)), \textup{point-\sot})$.
\item $f$ is continuous on $(K, \tau_{us*})$.
\item $f$ is continuous on $L_1$.
\end{enumerate}
\end{thm}

In classical Choquet theory, every state on $C(K)$ is given by integration against a Borel probability measure on $K$.
The restriction of a measure $\mu$ to $A(K)$ is a state on $A(K)$, and thus is evaluation at some point $x \in K$.
Then $x$ is called the barycentre of $\mu$ and $\mu$ is a representing measure for $x$.

In the nc setting, we replace states by nc states, namely \ucp\ maps $\mu:\rC(K) \to \B(H_n)$ for any cardinal $n$. 
As was observed above, the restriction of $\mu$ to $\rA(K)$ is a state on $\rA(K)$, and thus is evaluation at some $x\in K$, called the \emph{barycentre} of $\mu$. Likewise, we say that $\mu$ \emph{represents} $x$.
Stinespring's theorem \cite{St1955} says that a \ucp\ map $\mu$ has a dilation to a representation of $\rC(K)$, meaning that there is a $*$-representation $\pi$ of $\rC(K)$ on $H_m$ and an isometry $\beta \in M_{m,n}$ such that $\mu = \beta^* \pi \beta$. Moreover, there is a unique minimal dilation up to unitary equivalence fixing $H_n$.
Now there is a point $y \in K$ such that $\pi = \delta_y$; so we say that $(y,\beta)$ \emph{represents} $\mu$.
The uniqueness of the minimal representation $(y_0, \beta_0)$ for $y_0 \in H_{m_0}$ means that an arbitrary representation $(y,\beta)$ must decompose as $y = z_0 \oplus z_1$ with respect to $H_m = M \oplus M^\perp$ and $\beta = \begin{sbmatrix} \gamma_0 \\ 0 \end{sbmatrix}$ such that $\dim M = m_0$ and there is a unitary $\gamma \in \B(H_{m_0}, M)$ with $y_0 = \gamma z_0 \gamma^*$. We often write this as $y \simeq y_0 \oplus z_1$.
Note that this also shows that $\mu$ has a (unique) extension to a normal \ucp\ map on $\rB(K)$ as $\mu = \beta^* \delta_y \beta$.

In classical Choquet theory, the convex functions play a central role.
We require the analogue for nc functions.
Let $f  \in M_p(\rB(K))$ be self-adjoint bounded nc function. 
The \emph{epigraph} of $f$, $\epi(f)$, is defined by  
\[
 \epi_n(f) = \big\{ (x, \alpha) \in K_n \times M_p(M_n) : x \in K_n \AND \alpha \ge f(x) \big\} 
\]
as a subset of $\coprod_n K_n \times M_p(M_n)$.
Then $f$ is a \emph{convex nc function} if $\epi(f)$ is an nc convex set, and $f$ is \emph{lower semicontinuous} if $\epi(f)$ is closed.

\begin{prop} 
Let $K$ be a compact nc convex set and let $f$ be a self-adjoint bounded nc function in $M_p(\rB(K))$. 
The following are equivalent:
\begin{enumerate}
\item $f$ is convex.

\item \strut\qquad\qquad $f \big( \sum_i \alpha_i^* x_i \alpha_i \big) \le \sum_i \alpha_i^* f(x_i) \alpha_i$ \\
for $x_i \in K_{n_i}$ and $\alpha_i \in M_{n_i,n}$ such that $\sum_i \alpha_i^* \alpha_i = 1_n$.

\item \strut\qquad\qquad $f(\alpha^* x \alpha) \leq (1_p \otimes \alpha^*) f(x) (1_p \otimes \alpha) $\\ 
for every $m,n$, every $x \in K_m$ and every isometry $\alpha \in M_{m,n}$.

\item $f|_{K_n}$ is a scalar convex function for every $n$; i.e.,
\[
 f(t x + (1-t) y) \le t f(x) + (1-t) f(y) 
\]
for all $n$, all $x,y \in K_n$ and all $t \in [0,1]$.
\end{enumerate}
\end{prop}

The equivalence of the first three statements is straightforward, and these imply (4).
The equivalence of (4) is a generalization of the Hansen-Pedersen Jensen inequality \cite{HanPed2003}, 
which corresponds to the case of $K_1 = [0,1]$.

In the classical case, it is easy to see that $P(K) - P(K) = C_\bR(K)$.
And hence two measures which agree on the convex functions are equal.
We cannot show that every $f \in \rC(K)$ is the difference of two continuous nc convex functions, however we can get the desired consequence:

\begin{prop} 
Let $K$ be a compact nc convex set.
Suppose that $\mu, \nu : \rC(K) \to \M$ are \ucp\ maps such that $\mu(f) = \nu(f)$ for every $n$ and every convex nc function $f \in M_n(\rC(K))$. 
Then $\mu = \nu$.
\end{prop}

In the classical case, if $f,g$ are convex functions on a convex set $K$, then $\epi(f) \cap \epi(g) = \epi(f\vee g)$.
However in the non-commutative situation, since once cannot take the max of two self-adjoint operators, the same formula may not yield the epigraph of a function. For this reason, we use the notion of a multivalued (or set valued) function which takes care of this difficulty.

Let $K$ be an nc convex set and let $F : K \to \cM_n(\cM)_{sa}$ be a multivalued self-adjoint function. 
We say that $F$ is a \emph{multivalued nc function} if it is non-degenerate, graded, unitarily equivariant, preserves direct sums and is upwards directed, meaning that
\begin{enumerate}
\item $F(x) \ne \emptyset$ for every $x \in K$,
\item $F(K_m) \subseteq M_n(M_m)$ for all $m$,
\item $F(x \oplus y) = F(x) \oplus F(y)$ for every $x,y \in K$,
\item $F(\beta x \beta^*) = (1_n \otimes \beta) F(x) (1_n \otimes \beta^*)$ for every $x \in K_m$ and every unitary $\beta \in M_m$,
\item $F(x) = F(x) + M_n(M_m)_+$ for every $m$ and every $x \in K_m$.
\end{enumerate}

Say that $F$ is \emph{bounded} if there is a constant $\lambda > 0$ such that for every $\beta \in F(x)$, there is $\alpha \in F(x)$ with $\alpha \leq \beta$ such that $\|\alpha\| \leq \lambda$. If $F$ is bounded, then we let $\|F\|$ denote the infimum of all $\lambda$ that work. Otherwise we write $\|F\| = \infty$. If $G : K \to M_n(M)$ is another multivalued nc function, then we will write $F \leq G$ if $F(x) \supseteq G(x)$ for every $x \in K$.
Define
\[
\graph_m(F) = \{(x,\alpha) \in K_m \times \cM_n(\cM_m) : x \in K_m \text{ and } \alpha \in F(x) \}.
\]
We say that $F$ is \emph{convex} if $\graph(F) = \coprod_m \graph_m(F)$ is an nc convex set, and that $F$ is \emph{lower semicontinuous} if $\graph{F}$ is closed.
The \emph{convex envelope} of a bounded multivalued function $F : K \to M_n(\M)$ is the multivalued nc function $\env{F} : K \to M_n(\M)$ determined by the property
\[
 \graph{\env{F}}= \ol{\ncconv}(\graph(F)).
\]
That is, the graph of $\env{F}$ is the closed nc convex hull of the graph of $F$; whence $\env{F}$ is a convex lower semicontinuous multivalued nc function.

If $f$ is an nc function, then $\epi(f)$ is the graph of the multifunction $F(x) = [f(x), \infty)$; and $F$ is convex if and only if $f$ is convex.
We define $\env{f} := \env{\epi(f)}$.

It is straightforward to see that $\env{F} \le F$ and it is the largest convex lower semicontinuous multivalued nc function with this property.
Also if $F$ is bounded, then $\env{F}$ has the same bound.

In the classical case, the Hahn-Banach separation theorem shows that every lower semicontinuous function is the supremum of all smaller affine functions.
An analogue of this in the nc context is true, but is surprisingly difficult.

\begin{thm} \label{thm:convex-equals-upper-env}
Let $K$ be a compact nc convex set and let $F : K \to \cM_n(\cM)$ be a  bounded multivalued nc function. Then for $x \in K_p$, 
\[
\env{F}(x) = \bigcap_{m \in \bN\,} \bigcap_{a \leq 1_m \otimes F} \big\{ \alpha \in (\cM_n(\cM_p))_{sa} :  1_m \otimes \alpha \ge a(x) \big\},
\]
where the intersection is taken over all $m$ and all self-adjoint affine nc functions $a \in \cM_m(\cM_n(\rA(K)))_{sa}$ satisfying $a \leq 1_m \otimes F$. 
The same holds if we intersect over all $m \le \kappa$.
\end{thm}

\section{Orders on CP maps}

In classical Choquet theory, one puts a partial order on the probability measures on $K$ by setting $\mu \prec \nu$ if $\mu(f) \le \nu(f)$ for every continuous convex function $f$.
The existence of maximal measures in this ordering leads to a proof of the Choquet--Bishop--de Leeuw theorem, providing a representation for each point in $x$ as a measure which, in a technical sense, is supported on the set of extreme points \cite{BdL1959}.

In our setting, there are two natural partial orders on the \ucp\ maps (or nc states) from $\rC(K)$ to $H_n$ for any $n$.
Say that $\mu$ is dominated by $\nu$ in the \emph{nc Choquet order} and write $\mu \prec_c \nu$ if $\mu(f) \leq \nu(f)$ for every $n$ and every convex nc function $f \in \cM_n(\rC(K))$.
Say that $\mu$ is dominated by $\nu$ in the \emph{dilation order} and write $\mu \prec_d \nu$ if there are representations $(x,\alpha) \in K_n \times \cM_{n,m}$ for $\mu$ and $(y,\beta) \in K_p \times \cM_{p,m}$ for $\nu$ along with an isometry $\gamma \in \cM_{p,n}$ such that $\beta = \gamma \alpha$ and $x = \gamma^* y \gamma$.
We eventually show that these two orders are equal, and the interplay between them provides an additional tool in the noncommutative setting.

Note that these partial orders can only compare two nc states $\mu, \nu$ with the same range; i.e. $\mu, \nu \in L_n$ where $L$ is the nc state space of $\rC(K)$. If $a$ is nc affine, then $\pm a$ are nc convex; whence it follows that both orders only compare states with the same barycentre.

It is easy to see that the nc Choquet order behaves well under taking convex combinations and limits.
For each nc function $f$, we obtain that $\mu \prec_c \nu$ implies that $\mu(\env{f}) \le \nu(\env{f})$.
Since $\env{f}$ is a multifunction, this requires Theorem~\ref{thm:convex-equals-upper-env}.
A straightforward Zorn's lemma argument shows that for any $\mu$, there is $\mu \prec_c \nu$ such that $\nu$ is maximal in the nc Choquet order.

For the dilation order, there is a useful reformulation:

\begin{prop} \label {P:diln order}
Let $K$ be a compact nc convex set and let $\mu,\nu \in L_n$.
Then $\mu \prec_d \nu$ if and only if there is a representation $(x,\alpha) \in K_n \times M_{n,m}$ for $\mu$ 
and a \ucp\ map $\tau : \rC(K) \to M_n$ with barycentre $x$ satisfying $\nu = \alpha^* \tau \alpha$.
Moreover such a $\tau$ exists for any representation $(x,\alpha)$.
\end{prop}

The dilation order is convenient for computing \ucp\ maps of $\env{f}$.

\begin{thm} \label {T:convex-env-ucp-dilation-order}
Let $K$ be a compact nc convex set. 
Let $f \in M_n(\rB(K))$ be a bounded self-adjoint nc function with convex envelope $\env{f}$. 
Then for a \ucp\ map $\mu : \rC(K) \to M_k$,
\[
 \mu(\env{f}) = \bigcup_{\mu \prec_d \nu} [\nu(f), +\infty),
\]
where the union is taken over all $\nu$ such that $\mu \prec_d \nu$. 
Moreover, $\mu$ is maximal in the dilation order if and only if $ [\mu(f),+\infty) = \mu(\env{f})$ for all $f \in \rC(K)$.
\end{thm}

With these tools, we establish the equivalence of orders:

\begin{thm} \label {T:equivalence of orders}
Let $K$ be a compact nc convex set and let $\mu,\nu : \rC(K) \to M_n$ be \ucp\ maps. 
Then $\mu \prec_c \nu$ if and only if $\mu \prec_d \nu$.
\end{thm}

One interesting consequence of these ideas is that they lead to new, quite different and simpler, proofs of the Dritschel--McCullough and Davidson--Kennedy dilation theorems.

\section{Representation theorems}

The classical Choquet theorem states that if $K$ is a metrizable compact convex set, then every point in $K$ has a representing measure supported on $\partial K$. In the non-metrizable setting, there is the complication that the set of extreme points need not be Borel. Bishop and de Leeuw get around this by showing that there is a Baire measure $\mu$ representing $x \in K$ such that $\mu(f) = 0$ for every Baire function that vanishes on $\partial K$. Since the Baire functions form the $\sigma$-algebra of functions generated by $C(K)$, this turns out to be an appropriate solution.

In the non-commutative setting, Pedersen (see \cite{Pedersen}*{section 4.5}) defines what is known as the \emph{Pedersen--Baire envelope} of a C*-algebra $\fA$ as the monotone sequential completion $\rBI(\fA)$ of $\fA$ in $\fA^{**}$. We let $\rBI(K) = \rBI(\rC(K))$. This is the smallest C*-subalgebra of $\rB(K)$ containing $\rC(K)$ such that  the \sot\ limit of bounded monotone increasing sequences remain in $\rB^\infty(K)$. When $\fA$ is commutative, this algebra is the space of Baire functions. Every \ucp\ map on $\rC(K)$ has a unique sequentially normal extension to $\rBI(K)$.
We say that an nc state $\mu$ on $\rC(K)$ is \emph{supported on $\partial K$} if $\mu(f) = 0$ for all $f \in \rBI(K)$ such that $f|_{\partial K} = 0$.

The set $L_1$ of scalar states on $\rC(K)$ lie in the predual of $\rB(K)$; so these states are normal states on $\rB(K)$.
If $f \in \rB(K)$, the function $\hat f(\mu) = \mu(f)$ is bounded and affine on $L_1$.
Say that $f$ satisfies the \emph{barycentre formula} if whenever $\rho$ is a regular Borel probability measure on $L_1$ with barycentre $\mu$, then
\[
 \mu(f) = \int_{L_1} \nu(f) \,d\rho(\nu) .
\]
In general, this is not the case. It is valid however for $f \in \rBI(K)$.

We obtain the analogue of the classical representation theorem:

\begin{thm}[NC Choquet--Bishop--de Leeuw Theorem] \label {T:ncCBdL} \strut\\
Let $K$ be a compact nc convex set.
For each $x \in K_n$, there is a \ucp map $\mu:\rC(K) \to M_n$ with barycentre $x$ which is supported on $\partial K$.
\end{thm}

In the separable (metrizable) case of the classical theorem, Choquet produced a measure supported on the extreme points, which is a $G_\delta$ set.
In the separable case of the noncommutative theory, there is also an integral representation. It requires the machinery of integration against measures on the space of completely positive normal maps developed by Fujimoto \cite{Fuj1994}. We refer to \cite{DK_ncCh} for more detail.

Let $K$ be a compact nc convex set such that $\rA(K)$ is separable. For  $n \le \aleph_0$, a \emph{$\cM_n$-valued finite nc measure} on $K$ is a sequence 
$\lambda = (\lambda_m)_{m \le \aleph_0}$ such that each $\lambda_m$ is a $\ncpmaps(\cM_m,\cM_n)$-valued Borel measure and the sum
\[
\sum_{m \leq \aleph_0} \lambda_m(K_m)(1_m) \in \cM_n
\]
is weak* convergent. For $E \in \bor(K)$, we define $\lambda(E)$ by
\[
\lambda(E) = \sum_{m \leq \aleph_0} \lambda_m(E_m). 
\]
Say that $\lambda$ is \emph{supported on the extreme boundary $\partial K$} if 
\[
 \lambda_m(K_m \setminus \partial K) = 0_{\cM_m,\cM_n} \qforal m\le\aleph_0 .
\] 
Say that $\lambda$ is a \emph{$\cM_n$-valued nc probability measure} on $K$ if the above sum is equal to $1_n$. Finally, we will say that $\lambda$ is \emph{admissible} if each $\lambda_m$ is absolutely continuous with respect to a scalar-valued measure on $K_m$.

\begin{thm}[Noncommutative integral representation theorem] \label{thm:nc-choquet-theorem} \strut\\
Let $K$ be a compact nc convex set such that $\rA(K)$ is separable. Then for $x \in K$ there is an admissible nc probability measure $\lambda$ on $K$ that represents $x$ and is supported on $\partial K$, i.e. such that
\[
a(x) = \int_K a\, d\lambda, \qfor a \in \rA(K).
\]
\end{thm}

\section{Noncommutative simplices}

A compact convex set $K$ is a (Choquet) simplex if every point $x \in K$ has a unique representing measure that is maximal in the Choquet order. This motivates the following definition of a nc simplex, introduced by the second author and Shamovich in \cite{Kennedy-Shamovich}.

\begin{defn}
A compact nc convex set $K$ is an nc simplex if every point in $K$ has a unique representing map that is maximal in the nc Choquet order.
\end{defn}

Noncommutative simplices generalize classical simplices. Before stating this formally, we recall that for a classical compact convex set $C$, there is a unique compact nc convex set $K = \operatorname{max}(C)$ that is maximal with respect to the property that $K_1 = C$. The nc convex set $K$ is the nc state space of the operator system obtained by equipping the function system $\rA(C)$ of continuous affine functions on $C$ with its minimal operator system structure (cf. \cite{Paulsen-Todorov-Tomforde}). 

\begin{thm}
Let $K$ be a compact nc convex set. Then $K$ is an nc simplex with $\partial K = \partial K_1$ if and only if $K_1$ is a classical simplex, in which case $K = \operatorname{max}(K_1)$. Hence every classical simplex $C$ uniquely determines an nc simplex, namely $\operatorname{max}(C)$.
\end{thm}

\begin{cor} \label{cor:unique-operator-system-structure}
If $C$ is a classical simplex, then there is a unique operator system structure on the function system $\rA(C)$. 
\end{cor}

The converse of Corollary \ref{cor:unique-operator-system-structure} is known to hold for compact convex subsets of finite dimensional vector spaces (see \cite{Passer-Shalit-Solel}). However, it is an open problem whether it holds in general.

Many important properties of nc simplices can be derived from the following result.


\begin{thm}\label{thm:nc-simplex-map}
A compact nc convex set $K$ is an nc simplex if and only if there is a unital completely positive map $\phi : \rC(K) \to \rA(K)^{**}$ that restricts to the identity on $\rA(K)$, which is equivalent to the commutativity of the following diagram:
\[\begin{tikzcd}
	{\mathrm{C}(K)} && {\mathrm{A}(K)^{**}} \\
	& {\mathrm{A}(K)}
	\arrow["\varphi", from=1-1, to=1-3]
	\arrow[hook, from=2-2, to=1-1]
	\arrow[hook, from=2-2, to=1-3]
\end{tikzcd}\]
\end{thm}

Namioka and Phelps \cite{Namioka-Phelps} established a tensor product characterization of classical simplices. Specifically, they proved that a compact convex set $C$ is a simplex if and only if the function system $\rA(C)$ of continuous affine fucntions on $C$ is nuclear, meaning that
\[
\rA(C) \otimes_{\operatorname{min}} F = \rA(C) \otimes_{\operatorname{max}} F
\]
for every function system $F$, where $\otimes_{\operatorname{min}}$ and $\otimes_{\operatorname{min}}$ denote the minimal and maximal tensor product of functions systems respectively (see e.g. \cite{Paulsen-Todorov-Tomforde}). Effros \cite{Effros} proved that $C$ is a simplex if and only if the bidual $\rA(C)^{**}$ is injective in the category of function systems.

The results of Namioka-Phelps and Effros generalize to nc simplices.

\begin{thm} \label{thm:characterization-nc-simplices}
Let $K$ be a compact nc convex set. The following are equivalent:
\begin{enumerate}
\item $K$ is an nc simplex.
\item $\rA(K)^{**}$ is a von Neumann algebra.
\item $\rA(K)$ is $(\operatorname{c},\operatorname{max})$-nuclear, meaning that for every operator system $S$,
\[
\rA(K) \otimes_{\operatorname{c}} S = \rA(K) \otimes_{\operatorname{max}} S,
\]
where $\otimes_{\operatorname{c}}$ and $\otimes_{\operatorname{max}}$ denote the commuting and maximal tensor products of operator systems respectively (cf. \cite{KavPauTodTom}). 
\end{enumerate}
\end{thm}

There is an abundance of non-classical nc simplices, as the following examples demonstrate.

\begin{example}
Let $A$ be a unital C*-algebra with nc state space $K$, so that $\rA(K)$ is unitally completely order isomorphic to $A$. Then $\rA(K)^{**}$ is unitally completely order isomorphic to $A^{**}$, which is a von Neumann algebra. Hence by Theorem \ref{thm:characterization-nc-simplices}, $K$ is an nc simplex.
\end{example}

\begin{example}
Let $S$ be an operator system that is nuclear, or more generally has the weak expectation property (WEP). Then there is an injective operator system $T$ with $S \subseteq T \subseteq S^{**}$. Let $K$ denote the nc state space of $S$, so that $\rA(K)$ and $\rA(K)^{**}$ can be identified with $S$ and $S^{**}$. By the injectivity of $T$, the identity map on $\rA(K)$ extends to a unital completely positive map $\phi : \rC(K) \to T \subseteq \rA(K)^{**}$. Hence by Theorem \ref{thm:nc-simplex-map}, $K$ is an nc simplex.
\end{example}

\begin{example} 
Kirchberg and Wassermann \cite{KirWas1998} constructed a number of examples of operator systems that satisfy condition (2) in Theorem \ref{thm:characterization-nc-simplices}, which they called C*-systems. By Theorem \ref{thm:characterization-nc-simplices}, the nc state space of a C*-system is an nc simplex. They proved that a separable operator system $S$ is a C*-system if and only if there is a separable unital C*-algebra $A$ and a closed left ideal $L$ of $A$ such that $S$ is order isomorphic to the operator system $A/(L + L^*)$.
\end{example}

A classical simplex $C$ is said to be a Bauer simplex if its extreme boundary $\partial C$ is closed. The prototypical example of a Bauer simplex is the set $\operatorname{Prob}(X)$ of probability measures on a compact Hasudorff space equipped with the weak* topology, and Bauer \cite{Bauer} proved that every Bauer simplex arises in this way. Equivalently, since $\operatorname{Prob}(X)$ is the state space of the commutative C*-algebra $\rC(X)$ of continuous functions on $X$, a compact convex set $C$ is the state space of a unital commutative C*-algebra if and only if $C$ is a Bauer simplex.

A classical simplex $C$ is said to be a Poulsen simplex if its extreme boundary $\partial C$ is dense in $C$. Lindenstrauss, Olsen and Sternfeld \cite{Lindenstrauss-Olsen-Sternfeld} proved that there exists a unique metrizable Poulsen simplex.

There are good noncommutative generalizations of Bauer simplices and Poulsen simplices. In order to define these objects, it will be necessary to introduce a new topology on certain subsets of compact nc convex sets.

Let $K$ be a compact nc convex set. A point $x \in K$ is reducible if there are points $y,z \in K$ such that $x = y \oplus z$. Otherwise, $x$ is said to be irreducible. Let $\operatorname{Irr}(K)$ denote the set of irreducible points in $K$. Observe that $\partial K \subseteq \operatorname{Irr}(K)$. 

There is a bijective correspondence between points in $\operatorname{Irr}(K)$ and irreducible representations of $\rC(K)$ via the map $x \to \delta_x$. This correspondence induces a topology on $\operatorname{Irr}(K)$, namely the pullback of the hull-kernel topology on the irreducible representations of $\rC(K)$. We refer to this as the spectral topology on $\operatorname{Irr}(K)$.

\begin{defn}
A compact nc convex set $K$ is an nc Bauer simplex if it is an nc simplex and the extreme boundary $\partial K$ is closed in the set $\operatorname{Irr}(K)$ of irreducible points of $K$ with respect to the spectral topology.
\end{defn}

The next result characterizes the compact nc convex sets that arise as the nc state space of a unital C*-algebra.

\begin{thm}
A compact nc convex set is an nc Bauer simplex if and only if it is the nc state space of a unital C*-algebra.
\end{thm}

\begin{defn}
A compact nc convex set $K$ is an nc Poulsen simplex if it is an nc simplex and the extreme boundary $\partial K$ is dense in the set $\operatorname{Irr}(K)$ of irreducible points of $K$ with respect to the spectral topology. 
\end{defn}

\begin{thm}
A compact nc convex set $K$ is an nc Poulsen simplex if and only if $\rC^*_{\operatorname{max}}(\rA(K)) = \rC^*_{\operatorname{min}}(\rA(K))$.
\end{thm}

Lupini \cite{Lupini2} proved that there is a unique nuclear operator system $S$ with the property that $\rC^*_{\operatorname{max}}(\rA(K)) = \rC^*_{\operatorname{min}}(\rA(K))$. However, results of Kirchberg and Wassermann imply the existence of many non-nuclear operator systems with this property. In fact, their results imply that every compact nc convex set is the image of a nc Poulsen simplex under a continuous affine nc map (see \cite{Kennedy-Shamovich}).

In Section \ref{sec:dynamics}, we will see that nc simplices play an important role in noncommutative dynamics.

\begin{problem}
Let $K$ be an nc simplex such that $\rA(K)$ is separable. For $x \in K$, let $\lambda$ be an nc probability measure on $K$ that represents $x$ and is supported on the extreme boundary $\partial K$ as in Theorem \ref{thm:nc-choquet-theorem}. Is $\lambda$ unique in some sense?
\end{problem}

\section{Applications to approximation theory}

Classical Choquet theory applies to ordinary compact convex sets $K$, and the relationship between $A(K)$ and $C(K)$. The latter is a commutative C*-algebra. 
Kadison \cite{Kad1951} defined a \emph{concrete function system} to be a (closed) unital self-adjoint subspace of a unital commutative C*-algebra $C(X)$ which separates points of $X$.
He also defined function systems abstractly and showed that they are order isomorphic to a concrete function system.
A state on a function system is a unital positive linear functional.
If $A$ is a function system, the state space $K := S(A)$ of all states on $A$ is a compact, convex set in the weak-$*$ topology.
He shows that $A$ is order isomorphic to $A(K)$ via the map taking $a$ to $\hat{a}$, where $\hat{a}(x) = x(a)$.
Moreover, there is a $*$-homomorphism $q:C(K) \to C(X)$ such that $q(\hat a) = a$ for $a \in A$.
See \cite{AS2001}.

Note that $q^*:C(X)^* \to C(K)^*$ carries $S(A)$ onto $K$, considered as point evaluations.
In particular, point evaluations at points of $X$ are taken to points in $K$.
The Choquet boundary of $A(K)$ is just $\partial K$, the set of extreme points.
One then defines the Choquet boundary of $A$ to be $\partial_A X := (q^*)^{-1}(\partial K)$.

The classical Choquet order on finite positive measures on a compact convex set $K$ is defined as $\mu \prec \nu$ if $\mu(f) \le \nu(f)$ for all $f \in P(K)$, the space of continuous convex functions on $K$. 
Since whenever $a\in A(K)$, we have $\pm a \in P(K)$, we obtain $\mu(a) = \nu(a)$ whenever $\mu$ and $\nu$ are comparable.
In particular, two comparable probability measures have the same barycentre.
Bishop and de Leeuw \cite{BdL1959} were able to generalize Choquet's theorem to the non-separable case by establishing the existence of maximal measures dominating a probability measure $\mu$, and showing that it is supported (in a technical sense) on the set of extreme points.

In \cite{DK2021}, we define the dilation order in the commutative case.
We also introduce another order which yields new classical results.
In the following, given a finite positive measure $\mu$, the standard representation $\pi_\mu$ of $C(K)$ in $\B(L^2(\mu))$ is just $\pi_\mu(f) = M_f$, the multiplication operator. The range of $\pi_\mu$ is contained in a maximal abelian von Neumann algebra naturally isomorphic to $L^\infty(\mu)$.
The contant function $1_\mu$ lies in $L^2(\mu)$, and $\mu(f) = \ip{ \pi_\mu(f) 1_\mu, 1_\mu }$.  

\begin{defn} \label{defn:strong-dilation-order}
Let $K$ be a compact convex subset of a locally convex vector space. 
The {\em strong dilation order} ``$\prec_s$'' on $M^+(K)$ is defined for $\mu,\nu \in M^+(K)$ by $\mu \prec_s \nu$ if there is a positive map $\Phi:\rC(K) \to L^\infty(\mu)$ such that 
\[ \Phi(a) = \pi_\mu(a) \qforal a \in A(K) ,\]
and
\[ \nu(f) = \ip{\Phi(f) 1_\mu, 1_\mu} \qforal f \in \rC(K).\]
\end{defn}

A key result in \cite{DK2021} is the equivalence of the classical Choquet order and the strong dilation order.

\begin{thm}\label{T:choquet-equiv-strong-dilation-order}
Let $K$ be a compact convex subset of a locally convex vector space, and let $\mu,\nu \in M^+(K)$. 
Then $\mu \prec \nu$ if and only if $\mu \prec_s \nu$.
\end{thm}

A classical result of Korovkin states that if $\phi_n : C[0,1] \to C[0,1]$ is a sequence of positive maps such that $\lim_{n\to \infty} \| \phi_n(g) - g \| = 0$ for $g \in \{ 1,x,x^2\}$, then $\lim_{n\to \infty} \| \phi_n(f) - f \| = 0$ for all $f \in C[0,1]$.
\u Sa\u skin \cite{Sas1967} found a profound generalization of this.
If $X$ is compact and metrizable, a subset $G \subset C(X)$ is a \emph{Korovkin set} if whenever $\phi_n : C(X) \to C(X)$ is a sequence of positive linear maps satisfying $\lim_n \|\phi_n(g) - g\| = 0$ for all $g \in G$, then $\lim_n \|\phi_n(f) - f\| = 0$ for all $f \in C(X)$. 
The set $G$ determines a function system $A = \overline{\operatorname{span}(G \cup G^*)}$.
\u Sa\u skin established that a subset $G$ which separates points in $X$ is Korovkin if and only if $\partial_A X = X$.

One important consequence of the equivalence of the strong dilation order with the classical Choquet order is to remove the metrizability condition on $X$.

\begin{thm}\label{T:saskin}
Let $X$ be a compact Hausdorff space, and let $G$ be a subset of $\rC(X)$ containing the constant function $1$.
Let $A$ be the function system in $\rC(X)$ spanned by $G \cup G^*$. Then the following are equivalent:
\begin{enumerate}
\item $\partial_A X = X$.
\item $G$ is a Korovkin set.
\item For any compact Hausdorff space $Y$, unital $*$-homomorphism $\pi:\rC(X) \to \rC(Y)$ 
and net of positive linear maps $\phi_\lambda:\rC(X) \to \rC(Y)$ satisfying 
$\lim_\lambda \phi_\lambda(g) = \pi(g)$ for all $g \in G$, it follows that 
$\lim_\lambda \phi_\lambda(f) = \pi(f)$ for all $f \in \rC(X)$.
\end{enumerate}
\end{thm}

\section{Hyperrigidity}

Arveson \cite{Arv2011} established a stronger version of Korovkin's theorem.

\begin{thm}[Arveson 2011]
Let $\pi: C[0,1] \to \B(H)$ be a $*$-represent\-ation, and suppose that $\phi_n: C[0,1] \to \B(H)$ are unital positive linear maps such that
$\lim_{n\to\infty} \phi_n(g) = \pi(g)$ for $g \in \{ x, x^2\}$. Then $\lim_{n\to\infty} \phi_n(f) = \pi(f)$ for $f \in C[0,1]$.
\end{thm}

The important difference is that the range of the positive maps is $\B(H)$, not some commutative C*-algebra.
A positive map on a commutative C*-algebra is completely positive \cite{St1955}.
This led Arveson to make the following definition:

\begin{defn}
An operator system $S$ in an abstract C*-algebra $A = C^*(S)$ is \emph{hyperrigid} if for every (faithful) representation $\pi:A \to \B(H)$
and every sequence $\phi_n : A \to \B(H)$ of unital completely positive maps such that $\lim_{n\to\infty} \phi_n(s) = \pi(s)$ for $s \in S$,  one has
$\lim_{n\to\infty} \phi_n(a) = \pi(a)$ for $a \in A$.
\end{defn}

Many years earlier, Arveson \cite{Arv1969} had defined the \emph{unique extension property} for a representation $\pi$ of $A$ with respect to $S$
if there is a unique unital completely positive map $\phi:A \to \B(H)$ such that $\phi|_S = \pi|_S$.
In particular, and irreducible representation with the unique extension property was called a \emph{boundary representation}.
The set of all boundary representations plays the role of the Choquet boundary in the non-commutative setting.
While this idea was laid out in 1969, it wasn't completely realized until Arveson \cite{Arv2008} solved the separable case, and we in \cite{DK2015} established the existence of sufficiently many boundary representations in general.

Identifying hyperrigid operator systems is much easier because of the following:

\begin{thm}[Arveson 2011]
An operator system $S$ in a C*-algebra $A = C^*(S)$ is hyperrigid if and only if every representation of $A$ has the unique extension property relative to $S$.
\end{thm}

Arveson gave many explicit examples of hyperrigid operator systems, with various applications.
He formulated the following conjecture:

\begin{hyperconj}
Let $S$ be a separable operator system in a C*-algebra $A = C^*(S)$. Then $S$ is hyperrigid if and only if every irreducible representation of $A$ is a boundary representation for $S$.
\end{hyperconj}

Clearly, it is necessary that every irreducible representation be a boundary representation for $S$.

One can also see that in the commutative case, hyperrigidity is linked to the dilation order by Proposition~\ref{P:diln order}.
As strong dilation order requires a completely positive may into the von Neumann algebra $L^\infty(\mu)$ whereas dilation order only requires the range in $\B(H)$,
it follows that $\mu \prec \nu$ if and only if $\mu \prec_s \nu$ implies $\mu \prec_d \nu$.
It turns out that the converse fails.

In \cite{DK2021}, we establish

\begin{thm}
Let $K$ be a compact convex set and let $\mu$ be a probability measure on $K$.
Then the representation $\pi_\mu : C(K) \to \B(L^2(\mu))$ has the unique extension property relative to $A(K)$ if and only if $\mu$ is maximal on the dilation order.
\end{thm}

However we were not able to resolve the hyperrigidity conjecture in the commutative case. It comes down to the following:

\begin{problem}
If $A$ is a concrete function system generating a separable commutative C*-algebra $C(X)$ such that $\partial_A X = X$, is every measure which is Choquet maximal also maximal in the dilation order? 
\end{problem}

A positive answer confirms the conjecture in the commutative case, while a negative answer provides a counterexample.

It is now known that the hyperrigidity conjecture is false in general. 
Bilich and Dor-On \cite{Bilich-DorOn} found a counterexample for a C*-algebra which is an extension of the compacts by a commutative C*-algebra.
Then Bilich \cite{Bilich} analyzed the case of C*-correspondences and found a criterion which yields many other examples.

Clou\^{a}tre and Thompson \cite{Clouatre-Thompson} say that a representation $\pi$ of a C*-algebra $A$ has the \emph{unique tight extension property} 
with respect to $S$ if the only unital completely positive map $\phi: A \to \pi(A)''$ such that $\phi|_S = \pi|_S$ is $\pi$ itself.
This can be seen as being related to the strong dilation property. 
Note that when $\pi$ is irreducible, $\pi(A)'' = \B(H)$, so that the unique tight extension property and unique extension property coincide.
They prove

\begin{thm}[Clou\^{a}tre--Thompson]
Let $S$ be a separable operator system in a C*-algebra $A = C^*(S)$.
If every irreducible representation of $A$ is a boundary representation for $S$, then every representation of $A$ has the unique tight extension property.
\end{thm}

This leads to more results of the Korovkin--\u Sa\u skin type in the non-commutative setting generalizing Theorem~\ref{T:saskin}.
The commutative case of the hyperrigidity conjecture remains unresolved at this time.

\section{Dynamics} \label{sec:dynamics}

An action $G \curvearrowright X$ of a locally compact group $G$ on a compact Hausdorff space $X$ induces an affine action $G \curvearrowright \operatorname{Prob}(X)$ on the compact convex set of probability measures $\operatorname{Prob}(X)$ on $X$. This fundamental fact underlies much of topological dynamics and ergodic theory. The theory of nc convexity provides a noncommmutative analogue of this fact, as well as noncommutative analogues of some important applications.

Recall that from Corollary~\ref{C:categorical duality}, when $S$ is an operator system, 
$\Aut(S)$ denotes the group of unital order automorphisms of $S$; and when $K$ is a compact, convex set, $\Aut(K)$ denotes the group of affine homeomorphisms of $K$ onto itself.

\begin{defn}
\begin{enumerate}
\item An nc dynamical system is a triple $(S,G,\sigma)$ consisting of an operator system $S$, a locally compact group $G$ and a group homomorphism $\sigma : G \to \operatorname{Aut}(S)$ with the property that that for every $s \in S$, the map $G \to S : g \to \sigma_g(s)$ is continuous.
\item An affine nc dynamical system is a triple $(K,G,\kappa)$ consisting of a compact nc convex set $K$, a locally compact group $G$ and a group homomorphism $\kappa : G \to \operatorname{Aut}(K)$ with the property that for every $x \in K$, the map $G \to K : g \to \kappa_g(x)$ is continuous.
\end{enumerate}
\end{defn}

It follows from the dual equivalence between operator systems and compact nc convex sets that an nc dynamical system $(S,G,\sigma)$ correpsonds to an affine nc dynamical system $(K,G,\kappa)$, where $K$ denotes the nc state space of $S$. The relationship between $\sigma$ and $\kappa$ is given by
\[
(\kappa_g(\phi))(s) = \phi(\sigma_{g^{-1}}(s)), \quad \text{for} \quad s \in S,\ g \in G, \phi \in K.
\]

Let $(A,G,\sigma)$ be a C*-dyamical system. This is a special case of an nc dyamical system where the operator system is a unital C*-algebra $A$.

If $A$ is commutative, say $A = \rC(X)$ for a compact Hausdorff space $X$, then $(\rC(X),G,\sigma)$ corresponds to an action $G \curvearrowright X$ and so gives rise to an induced action $G \curvearrowright \operatorname{Prob}(X)$ on the state space $\operatorname{Prob}(X)$ of $\rC(X)$. An important consequence of the fact that $\operatorname{Prob}(X)$ is a (Bauer) simplex is the fact that the compact convex subset of invariant measures
\[
\operatorname{Prob}(X)^K = \{\mu \in \operatorname{Prob}(X) : g \mu = \mu \text{ for all } g \in G\}
\]
is a simplex. The extreme points of $\operatorname{Prob}(X)^K$ are the ergodic measures, and the Choquet-Bishop-de Leeuw theorem combined with the fact that $\operatorname{Prob}(X)$ is a simplex implies the ergodic decomposition theorem: every element in $\operatorname{Prob}(X)$ has a unique representing measure supported on the ergodic measures.

If $A$ is noncommutative, then there is still an induced affine action on the state space of $A$. However, this action is typically no longer as useful. For example, one major issue is that the state space of a noncommutative C*-algebra is never a simplex. As a result, the set of invariant states is generally not a simplex.

Letting $K$ denote the nc state space of $A$, it is natural to instead consider affine nc dynamical system $(K,G,\sigma)$ corresponding to $(A,G,\sigma)$. By Theorem \ref{thm:characterization-nc-simplices}, $K$ is an nc simplex, so the following result implies that the compact nc convex subset
\[
K^G = \{x \in K : gx = x \text{ for all } g \in G \}
\]
of invariant nc states is an nc simplex.

\begin{thm}
Let $(K,G,\kappa)$ be an affine nc dynamical system. If $K$ is an nc simplex, the the compact nc convex subset $K^G$ of invariant nc states is an nc simplex.
\end{thm}

Glasner and Weiss proved a dynamical characterization of property (T). Specifically, they proved that a second countable locally compact group $G$ has property (T) if and only if for every action $G \curvearrowright X$ of $G$ on a compact Hausdorff space $X$, the compact convex subset $\operatorname{Prob}(X)^G$ of invariant probability measures is a Bauer simplex. By Bauer's characterization of state spaces of commutative C*-algebras, this is equivalent to the statement that $G$ has property (T) if and only if for every commutative C*-dynamical system $(\rC(X),G,\sigma)$, the compact convex subset of invariant states is the state space of a commutative C*-algebra. The following theorem is a noncommutative generalization of this result. 

\begin{thm}
Let $G$ be a second countable locally compact group. Then $G$ has Kazhdan's property (T) if and only if whenever $(A,G,\sigma)$ is a C*-dynamical system with corresponding affine nc dynamical system $(K,G,\kappa)$, the compact nc convex subset $K^G$ of invariant nc states is the nc state space of a unital C*-algebra.
\end{thm}

The ideas in this section have also played an important role in the study of the ideal structure of reduced crossed product C*-algebras. See for example \cite{Kennedy-Schafhauser}.


\end{document}